# ON SAMPLING OF STATIONARY INCREMENT PROCESSES

By J. M. P. Albin[1]

*Chalmers University of Technology*

Under a complex technical condition, similar to such used in extreme value theory, we find the rate $q(\varepsilon)^{-1}$ at which a stochastic process with stationary increments $\xi$ should be sampled, for the sampled process $\xi(\lfloor\cdot/q(\varepsilon)\rfloor q(\varepsilon))$ to deviate from $\xi$ by at most $\varepsilon$, with a given probability, asymptotically as $\varepsilon\downarrow 0$. The canonical application is to discretization errors in computer simulation of stochastic processes.

**1. Introduction.** Let $\{\xi(t)\}_{t\in[0,1]}$ be a stationary increment (SI) process. In Section 4, under Condition $\mathfrak{S}$ from Section 2, we prove that, for suitable functions $q_1, q_2, w, \overline{F} > 0$, a constant $\kappa > 0$ and an open interval $J \subseteq \mathbb{R}$, denoting $\lfloor x\rfloor \equiv \sup\{k\in\mathbb{Z}: k\le x\}$,

$$\lim_{\varepsilon\downarrow 0}\mathbf{P}\left\{\frac{\sup_{t\in[0,1]}[\xi(t)-\xi(\lfloor t/q_1(\varepsilon)\rfloor q_1(\varepsilon))]-\varepsilon}{w(\varepsilon)}\le x\right\} = e^{-\kappa\overline{F}(x)} \quad \text{for } x\in J. \tag{1.1}$$

$$\lim_{\varepsilon\downarrow 0}\mathbf{P}\left\{\frac{\sup_{t\in[0,1]}|\xi(t)-\xi(\lfloor t/q_2(\varepsilon)\rfloor q_2(\varepsilon))|-\varepsilon}{w(\varepsilon)}\le x\right\} = e^{-\kappa\overline{F}(x)} \quad \text{for } x\in J. \tag{1.2}$$

Methods of proof relate to Albin (1998), for the analysis of single sampling intervals, and to Albin (1990), when assembling these to the limits (1.1) and (1.2). The main inspirations for those works in turn were Berman (1982, 1992) and Pickands (1969a, b) together with Leadbetter, Lindgren and Rootzén [(1983), Chapters 12 and 13], respectively.

In Section 3 we discuss Condition $\mathfrak{S}$ and the role of (1.1) and (1.2) in applications.

---

Received October 2002; revised May 2003.
[1]Supported by NFR Grant M 5105-20005196/2000.

Versions of (1.1) and (1.2) have a history for Gaussian processes; see Belyaev and Simonyan (1979), Piterbarg and Seleznjev (1994), Seleznjev (1996) and Hüsler (1999). Section 5 features applications to stable processes. Future applications could be to light-tailed Lévy processes and Ornstein–Uhlenbeck processes [cf. Albin (2003a, b)].

Albeit similar in appearance, (1.2) has little if any relation to the small ball problem, to estimate $\mathbf{P}\{\sup_{t\in[0,1]}|\xi(t)|\le\varepsilon\}$ as $\varepsilon\downarrow 0$ [e.g., Kuelbs and Li (1993)], as (1.1) is unrelated to estimates of $\mathbf{P}\{\sup_{t\in[0,1]}\xi(t)\le\varepsilon\}$ [e.g., Marcus (2000)].

**2. Condition $\mathfrak{S}$.** The conditions (2.1)–(2.8) on the process $\{\xi(t)\}_{t\in[0,1]}$ listed below will be collectively referred to as Condition $\mathfrak{S}(i)$, with $i=1$ or $i=2$:

(a) Let $q=q(\varepsilon),\tilde q=\tilde q(\varepsilon),\hat q=\hat q(\varepsilon)$ and $w=w(\varepsilon)$ be strictly positive functions on $(0,\infty)$, with $q(\varepsilon),w(\varepsilon)\downarrow 0$ as $\varepsilon\downarrow 0$, and with either $\tilde q=1$ or $Q_1\equiv\lim_{\varepsilon\downarrow 0}1/\tilde q=\infty$, and either $\hat q=\tilde q$ or $Q_2\equiv\lim_{\varepsilon\downarrow 0}\tilde q/\hat q=\infty$. Set $I=(0,Q_1)$ or $I=(0,Q_1]$.

(b) Let $\xi$ be SI, separable and measurable, defined on a complete probability space. Depending on whether (1.1) or (1.2) is in focus, put $\xi_\varepsilon=\xi_{1,\varepsilon}$ or $\xi_\varepsilon=\xi_{2,\varepsilon}$, where

$$\begin{aligned}\xi_{1,\varepsilon}(t)&\equiv[\xi(qt)-\xi(q\lfloor t-0\rfloor)-\varepsilon]/w\\ \xi_{2,\varepsilon}(t)&\equiv[|\xi(qt)-\xi(q\lfloor t-0\rfloor)|-\varepsilon]/w\end{aligned}\qquad\text{for } t>0.$$

(c) Pick an open set $J\subseteq\mathbb{R}$ with $0\in J$, a continuous function $\overline{F}:J\to(0,\infty)$ with $\overline{F}(0)=1$, and an integrable function $\overline{G}:[0,Q_1)\to(0,\infty)$ with $\overline{G}(0)=1$. Assume that, for $0<s/Q_2<Q_1-r\le Q_1$ and $x\in J$,

$$\lim_{\varepsilon\downarrow 0}\frac{\mathbf{P}\{\xi_\varepsilon(1-r\tilde q-s\hat q)>x\}}{(\hat q/\tilde q)q}=\overline{F}(x)\overline{G}\biggl(r+\frac{s}{Q_2}\biggr).\tag{2.1}$$

(d) Assume that, for $T\in I$ and $x\in J$,

$$\begin{aligned}&\lim_{\varepsilon\downarrow 0}\int_0^{T\wedge(1/\tilde q)}\frac{\mathbf{P}\{\xi_\varepsilon(1-r\tilde q)>x\}}{\mathbf{P}\{\xi_\varepsilon(1)>0\}}\,dr\\ &\qquad=\overline{F}(x)\int_0^T\overline{G}(r)\,dr\equiv\overline{F}(x)E(T;\overline{G})>0.\end{aligned}\tag{2.2}$$

(e) Assume that, for $n\in\mathbb{N}$, $r\in(0,Q_1)$, $t\in(0,rQ_2)^n$ and (independently of) $x\in J$,

$$\lim_{\varepsilon\downarrow 0}\mathbf{P}\Biggl\{\bigcap_{i=1}^n\{\xi_\varepsilon(1-r\tilde q+t_i\hat q)>x\}\Bigg|\xi_\varepsilon(1-r\tilde q)>x\Biggr\}=P_r(t).\tag{2.3}$$

(f) If $Q_2 = \infty$, assume that, for any $x \in J$, a.e. for $r \in I$,

$$\lim_{\widehat{T}\to\infty} \limsup_{\varepsilon\downarrow 0} \int_{\widehat{T}}^{r\tilde{q}/\hat{q}} \mathbf{P}\{\xi_\varepsilon(1-\tilde{q}r+\hat{q}s) > x | \xi_\varepsilon(1-\tilde{q}r) > x\}\, ds = 0. \tag{2.4}$$

(g) Assume that, for $T \in I$ and $\hat{x} < x$ in $J$, with $\oint$ denoting the direct upper Riemann integral of nonnegative functions [e.g., Asmussen (1987)],

$$\begin{aligned} &\lim_{a\downarrow 0} \limsup_{\varepsilon\downarrow 0} \frac{1}{a\mathbf{P}\{\xi_\varepsilon(1) > 0\}} \\ &\qquad \times \oint_0^{T\wedge(1/\tilde{q})} \mathbf{P}\Big\{ \sup_{s\in[0,1]} \xi_\varepsilon((1-t\tilde{q}-sa\hat{q})^+) > x, \\ &\qquad\qquad \xi_\varepsilon(1-t\tilde{q}) \vee \xi_\varepsilon((1-t\tilde{q}-a\hat{q})^+) \le \hat{x} \Big\}\, dt = 0. \end{aligned} \tag{2.5}$$

(h) If $I = (0, Q_1)$, assume that, for $T \in I$ and $x \in J$,

$$\inf_{T\in I} \limsup_{\varepsilon\downarrow 0} \frac{\mathbf{P}\{\sup_{r\in[T\wedge(1/\tilde{q}),(1/\tilde{q}))} \xi_\varepsilon(1-r\tilde{q}) > x\}}{(\tilde{q}/\hat{q})\mathbf{P}\{\xi_\varepsilon(1) > 0\}} = 0. \tag{2.6}$$

(i) Assume that, for $T \in I$ and $x \in J$,

$$\lim_{h\downarrow 0} \limsup_{\varepsilon\downarrow 0} \sup_{r,s\in[0,T\wedge(1/\tilde{q}))} \sum_{n=2}^{\lfloor h/q \rfloor} \frac{q\mathbf{P}\{\xi_\varepsilon(1-r\tilde{q}) > x, \xi_\varepsilon(n-s\tilde{q}) > x\}}{\mathbf{P}\{\xi_\varepsilon(1-r\tilde{q}) > 0\}\mathbf{P}\{\xi_\varepsilon(1-s\tilde{q}) > 0\}} = 0. \tag{2.7}$$

(j) Given $T \in I$, $x \in J$ and $a, h > 0$, assume that, for any $s_1, \ldots, s_i, t_1, \ldots, t_{i'} \in [0, 1/q] \cap \bigcup_{n=1}^\infty (n - T\tilde{q}, n]$ with $s_{j+1} - s_j \ge a\hat{q}$, $t_{j+1} - t_j \ge a\hat{q}$ and $t_1 - s_i \ge h/q$,

$$\begin{aligned} &\lim_{\varepsilon\downarrow 0} \mathbf{P}\Big\{ \bigcap_{k=1}^{i} \{\xi_\varepsilon(s_k) \le x\}, \bigcap_{\ell=1}^{i'} \{\xi_\varepsilon(t_\ell) \le x\} \Big\} \\ &\qquad - \mathbf{P}\Big\{ \bigcap_{k=1}^{i} \{\xi_\varepsilon(s_k) \le x\} \Big\} \mathbf{P}\Big\{ \bigcap_{\ell=1}^{i'} \{\xi_\varepsilon(t_\ell) \le x\} \Big\} = 0. \end{aligned} \tag{2.8}$$

THEOREM 1. *If the Condition $\mathfrak{S}(i)$ holds, for an $i \in \{1, 2\}$, then $(1.i)$ holds for some constant $\kappa > 0$.*

**3. Discussion.** The idea behind (1.1) and (1.2) is to use extreme value technology, for the extreme values $\{\varepsilon + xw\}_{x\in J}$ of the processes $\xi(qt) - \xi(q\lfloor t-0 \rfloor)$ and $|\xi(qt) - \xi(q\lfloor t-0 \rfloor)|$, where $q = q(\varepsilon)$ is chosen small enough to make $\varepsilon$ an "extremely large" value for these processes as $\varepsilon \downarrow 0$. Conditions (2.1)–(2.8) emerge from inserting $\xi(qt) - \xi(q\lfloor t-0 \rfloor)$ and $|\xi(qt) - \xi(q\lfloor t - 0\rfloor)|$ in analogous conditions from extremes [e.g., Albin (1990, 1998)]. What makes life difficult here is that these processes depend on the level $\varepsilon$ and do not possess any invariances like self-similarity (SS) or SI.

However, for $\xi$ SSSI, (1.1) and (1.2) are easily transformed to more conventional extreme value problems, involving processes that do not depend on the level.

Schemes from extremes are natural take-off points to check Condition $\mathfrak{S}$. Basically, we refer readers to the literature on extremes for that purpose (but see Section 5). In fact, very much can be written on how to check (2.1)–(2.8) [see, e.g., Albin (1998) for a hint on this]. We plan to return to this in the future.

We discuss Condition $\mathfrak{S}(i)$ for $i=1$; the case when $i=2$ is very similar.

Condition (2.1) implies that

$$\mathbf{P}\{\xi_\varepsilon(1)>0\}=\mathbf{P}\{\xi(q)-\xi(0)>\varepsilon\}\sim(\hat{q}/\tilde{q})q \qquad \text{as } \varepsilon\downarrow 0 \tag{3.1}$$

and that

$$\begin{aligned}\lim_{\varepsilon\downarrow 0}\frac{\mathbf{P}\{\xi_\varepsilon(1)>x\}}{\mathbf{P}\{\xi_\varepsilon(1)>0\}}&=\mathbf{P}\left\{\frac{\xi(q)-\xi(0)-\varepsilon}{w}>x\,\middle|\,\frac{\xi(q)-\xi(0)-\varepsilon}{w}>0\right\}\\&=\overline{F}(x).\end{aligned}\tag{3.2}$$

This means that $\xi(q)-\xi(0)$ belongs to a domain of attraction of extremes [e.g., Albin (1990), Section 1], modulo the fact that this random variable depends on $\varepsilon$. In addition to (3.1) and (3.2), (2.1) requires that $\tilde{q}$ is chosen in such a way that

$$\lim_{\varepsilon\downarrow 0}\frac{\mathbf{P}\{\xi_\varepsilon(1-\tilde{q}r)>x\}}{\mathbf{P}\{\xi_\varepsilon(1)>0\}}=\overline{F}(x)\overline{G}(r) \qquad \text{where } \overline{G}\in\mathbb{L}^1. \tag{3.3}$$

Arguably, in the presence of (3.2), this is a reasonable requirement.

In the presence of (2.1), condition (2.2) is just a technicality. Notice that, by (2.1) and Fatou's lemma, (2.2) for a certain $T=\widehat{T}$ implies (2.2) for $T\le\widehat{T}$.

Typically, (2.3) is verified by means of establishing that

$$\{(\xi_\varepsilon(1-r\tilde{q}+t\hat{q})|\xi_\varepsilon(1-r\tilde{q})>x)\}_{t\ge 0}\xrightarrow{d}\{\zeta_x(t)\}_{t\ge 0} \qquad \text{as } \varepsilon\downarrow 0, \tag{3.4}$$

for some stochastic process $\{\zeta_x(t)\}_{t>0}$. This together with (2.1) implies (2.3), since

$$\begin{aligned}&\limsup_{\varepsilon\downarrow 0}\mathbf{P}\{\xi_\varepsilon(1-r\tilde{q}+t\hat{q})=x|\xi_\varepsilon(1-r\tilde{q})>x\}\\&\quad\le\limsup_{\delta\downarrow 0}\limsup_{\varepsilon\downarrow 0}\frac{\mathbf{P}\{\xi_\varepsilon(1-r\tilde{q}+t\hat{q})\in(x-\delta,x+\delta]\}}{\mathbf{P}\{\xi_\varepsilon(1-r\tilde{q})>x\}}\\&\quad=\limsup_{\delta\downarrow 0}[\overline{F}(x-\delta)-\overline{F}(x+\delta)]\left[\overline{G}\left(r+\frac{t}{Q_2}\right)\right]=0.\end{aligned}$$

In particular, (2.3) implies that

$$\lim_{\varepsilon\downarrow 0}\mathbf{P}\{\xi_\varepsilon(1-r\tilde{q}+t\hat{q})>x|\xi_\varepsilon(1-r\tilde{q})>x\}=P_r(t). \tag{3.5}$$

Assuming that (2.3) is established via (3.4), a formal argument gives

$$\lim_{\varepsilon\downarrow 0}\int_{\widehat{T}}^{r\tilde{q}/\hat{q}}\mathbf{P}\{\xi_\varepsilon(1-\tilde{q}r+\hat{q}s)>x\mid\xi_\varepsilon(1-\tilde{q}r)>x\}\,ds$$
$$=\int_{\widehat{T}}^{\infty}\mathbf{P}\{\zeta_x(s)>x\}\,ds=\mathbf{E}\biggl\{\int_{\widehat{T}}^{\infty}\mathbf{1}_{(x,\infty)}(\zeta_x(s))\,ds\biggr\}$$

when $Q_2=\infty$: the interpretation of (2.4) is thus one about asymptotically small mean sojourn times above the level $x$, for the process $\zeta_x$ at infinity.

Condition (2.5) is typically verified by bounds on increments of $\xi_\varepsilon$ (see, e.g., Section 5). For example, (2.5) holds if, for some constants $K,d>0$ and $e>1$,

(3.6) $$\mathbf{P}\{\xi_\varepsilon(1-r\tilde{q}+s\hat{q})\le x-\delta,\xi_\varepsilon(1-r\tilde{q})>x\}\le K\delta^{-d}s^e\mathbf{P}\{\xi_\varepsilon(1)>0\}$$

for $\varepsilon$ large enough, with uniformity in $r$. When (2.1) and (3.4) holds, we have

$$\lim_{\varepsilon\downarrow 0}\frac{\mathbf{P}\{\xi_\varepsilon(1-r\tilde{q}+s\hat{q})\le x-\delta,\xi_\varepsilon(1-r\tilde{q})>x\}}{\mathbf{P}\{\xi_\varepsilon(1)>0\}}=\overline{F}(x)\overline{G}(r)\mathbf{P}\{\zeta_x(s)\le x-\delta\}.$$

For $\overline{F}$ and $\overline{G}$ "nice," the interpretation of (3.6) thus is [recall that $\zeta_x(0)>x$]

$$\mathbf{P}\{\zeta_x(s)\le x-\delta\}\le K\delta^{-d}s^e.$$

Condition (2.5) is often easier to verify for $I=(0,Q_1)$ than for $I=(0,Q_1]$. The price for this is that (2.6) is required, to take care of $\xi_\varepsilon(t)$ for $t$ close to zero.

Condition (2.7) is a mixing conditon. It requires in essence that

$$\mathbf{P}\{\xi_\varepsilon(1-r\tilde{q})>x,\xi_\varepsilon(n-s\tilde{q})>x\}\le K\mathbf{P}\{\xi_\varepsilon(1-r\tilde{q})>0\}\mathbf{P}\{\xi_\varepsilon(n-s\tilde{q})>0\}$$

for some constant $K>0$. Condition (2.8) is also a mixing conditon that is imposed for process values that are farther apart than those in (2.7).

In most applications, $Q_2=1$ so that $\hat{q}=\tilde{q}$, since $Q_2=\infty$ usually corresponds to rather pathological processes: for example, for fractional Brownian motion, one has $Q_2=\infty$ for the index of self-similarity $H\in(0,\frac{1}{2})$, while $Q_2=1$ for $H\in[\frac{1}{2},1)$ [see Piterbarg and Seleznjev (1994) and Seleznjev (1996)].

For $Q_2=1$, the functions $q$, $w$ and $\hat{q}=\tilde{q}$ are typically determined from (3.1), (3.2) and (3.3), respectively (in that order). For $Q_2=\infty$ four functions are sought, and an additional equation is required, for example, (3.5). Of course, a rigorous use of Theorem 1 requires a verification of the entire Condition $\mathfrak{S}$.

The following two examples are canonical for the role of (1.1) in applications:

Example 1. In a computer simulation a largest deviation $d=\varepsilon+xw(\varepsilon)$ between the discretized simulated process $\xi(\lfloor\cdot/q(\varepsilon)\rfloor q(\varepsilon))$ and the process $\{\xi(t)\}_{t\in[0,1]}$ itself is tolerated, with a certain probability $p=e^{-\kappa\overline{F}(x)}$. For this, the sampling rate should be $q(\varepsilon)^{-1}$, asymptotically as $\varepsilon\downarrow 0$ (same thing as $d\downarrow 0$). Here the deviation is measured "two-sided" in (1.2), and "one-sided" in (1.1).

An applied risk application can be to estimate a high quantile

$$\mathbf{P}\Big\{\sup_{t\in[0,1]}\xi(t)>u\Big\}\leq 2p \qquad \text{for } u \text{ large and } p \text{ small.} \tag{3.7}$$

This is solved using computer simulations, to find a $u$ such that

$$\mathbf{P}\Bigg\{\bigcup_{k=0}^{\lfloor 1/q(\varepsilon)\rfloor}\{\xi(kq)>u-d\}\Bigg\}=p=1-e^{-\kappa\overline{F}(x)} \qquad \text{where } d=\varepsilon+xw(\varepsilon),$$

for $\varepsilon$ small enough that the limit (1.1) is accurate, because then (1.1) gives

$$\begin{aligned}\mathbf{P}\Big\{&\sup_{t\in[0,1]}\xi(t)>u\Big\}\\ &\leq\mathbf{P}\Bigg\{\bigcup_{k=0}^{\lfloor 1/q(\varepsilon)\rfloor}\{\xi(kq)>u-d\}\Bigg\}\\ &\quad+\mathbf{P}\Big\{\sup_{t\in[0,1]}[\xi(t)-\xi(\lfloor t/q(\varepsilon)\rfloor q(\varepsilon))]>\varepsilon+xw(\varepsilon)\Big\}\\ &\leq p+1-e^{-\kappa\overline{F}(x)}=2p.\end{aligned}$$

Example 2. In a theoretical risk application, the problem is again to solve (3.7). Let $\xi$ be stationary, with $\xi(0)$ belonging to a domain of attraction, that is,

$$\lim_{u\to\infty}\frac{\mathbf{P}\{\xi(0)>u+y\tilde{w}(u)\}}{\mathbf{P}\{\xi(0)>u\}}=\overline{H}(y) \qquad \text{for } y\in J, \tag{3.8}$$

for some functions $\tilde{w},\overline{H}>0$. Picking a $-y\in J\cap(-\infty,0)$, (3.7) is solved taking

$$y\tilde{w}(u)=\varepsilon+xw(\varepsilon) \quad \text{and} \quad \Big(\frac{1}{q(\varepsilon)}+1\Big)H(-y)\mathbf{P}\{\xi(0)>u\}=p=1-e^{-\kappa\overline{F}(x)},$$

for $\varepsilon$ small enough and $u$ large enough that the limits (1.1) and (3.8) are accurate, because then we have, for such $\varepsilon$ and $u$, by stationarity, (1.1) and (3.8),

$$\mathbf{P}\Big\{\sup_{t\in[0,1]}\xi(t)>u\Big\}$$

$$
\begin{aligned}
&\le \mathbf{P}\Biggl\{\bigcup_{k=0}^{\lfloor 1/q(\varepsilon)\rfloor}\{\xi(kq) > u - y\tilde{w}(u)\}\Biggr\} \\
&\quad + \mathbf{P}\Biggl\{\sup_{t\in[0,1]}\biggl[\xi(t) - \xi\biggl(\biggl\lfloor\frac{t}{q(\varepsilon)}\biggr\rfloor q(\varepsilon)\biggr)\biggr] > \varepsilon + xw(\varepsilon)\Biggr\} \\
&\le \biggl(\frac{1}{q(\varepsilon)} + 1\biggr)\mathbf{P}\{\xi(0) > u - y\tilde{w}(u)\} + 1 - e^{-\kappa\overline{F}(x)} = 2p.
\end{aligned}
$$

We conclude this section with an example featuring, possibly, the only nontrivial process for which the probability in (1.1) can be calculated explicitly.

EXAMPLE 3. Let $B$ be Brownian motion (arbitrarily started), $W \equiv B - B(0)$ and

$$
q_i(\varepsilon) \equiv \frac{\varepsilon^2}{\operatorname{Var}\{W(1)\}}\biggl[4\ln\biggl(\frac{1}{\varepsilon}\biggr) + \ln\biggl(4\ln\biggl(\frac{1}{\varepsilon}\biggr)\biggr) + 2\ln\biggl(\frac{2\operatorname{Var}\{W(1)\}i}{\sqrt{2\pi}}\biggr)\biggr]^{-1}.
$$

By self-similarity and other elementary properties of Brownian motion, we get

$$
\begin{aligned}
&\mathbf{P}\Biggl\{\sup_{t\in[0,1]} B(t) - B\biggl(\biggl\lfloor\frac{t}{q}\biggr\rfloor q\biggr) \le \varepsilon\Biggr\} \\
&\quad = \biggl(\mathbf{P}\Bigl\{\sup_{t\in[0,q]} W(t) \le \varepsilon\Bigr\}\biggr)^{\lfloor 1/q\rfloor}\mathbf{P}\Bigl\{\sup_{t\in[0,1-\lfloor 1/q\rfloor q]} W(t) \le \varepsilon\Bigr\} \\
&\quad = (1 - 2\mathbf{P}\{W(q) > \varepsilon\})^{\lfloor 1/q\rfloor}\biggl(1 - 2\mathbf{P}\biggl\{W\biggl(1 - \biggl\lfloor\frac{1}{q}\biggr\rfloor q\biggr) > \varepsilon\biggr\}\biggr) \\
&\quad \sim \biggl(1 - \frac{\sqrt{2Cq}}{\sqrt{\pi}\varepsilon}\exp\biggl\{-\frac{\varepsilon^2}{2Cq}\biggr\}(1 + o(1))\biggr)^{1/q} \to e^{-2}.
\end{aligned}
$$

**4. Proof of Theorem 1.** The sojourn time

$$
L_\varepsilon(s;x) \equiv \int_0^{s\wedge(1/\tilde{q})}\mathbf{1}_{(x,\infty)}(\xi_\varepsilon(1-\tilde{q}r))\,dr, \qquad (s,x) \in I \times J,
$$

satisfies

$$
\mathbf{E}\{L_\varepsilon(T;x)\} \sim \mathbf{P}\{\xi_\varepsilon(1) > 0\}\overline{F}(x)E(T;\overline{G}) \qquad \text{as } \varepsilon \downarrow 0 \tag{4.1}
$$

for $T \in I$ and $x \in J$ [by (2.2) and Fubini's theorem]. By the elementary identity

$$
\begin{aligned}
&\mathbf{1}_{(y,\infty)}(L_\varepsilon(T;x))\int_0^{T\wedge(1/\tilde{q})}\mathbf{1}_{(-\infty,y]}(L_\varepsilon(s;x))\mathbf{1}_{(x,\infty)}(\xi_\varepsilon(1-\tilde{q}s))\,ds \\
&\quad = \mathbf{1}_{(y,\infty)}(L_\varepsilon(T;x))y
\end{aligned}
$$

for $T \in I$ and $y > 0$, it therefore follows that

$$\frac{\hat{q}}{\tilde{q}} \int_y^\infty \mathbf{P}\Big\{\tilde{q} L_\varepsilon(T;x)/\hat{q} > z\Big\}\, dz$$

$$= \int_0^\infty \mathbf{P}\{(L_\varepsilon(T;x) - \hat{q}y/\tilde{q}) > \hat{z}\}\, d\hat{z}$$

$$= \mathbf{E}\{(L_\varepsilon(T;x) - \hat{q}y/\tilde{q})\mathbf{1}_{(\hat{q}y/\tilde{q},\infty)}(L_\varepsilon(T;x))\} \tag{4.2}$$

$$= \mathbf{E}\Big\{\mathbf{1}_{(\hat{q}y/\tilde{q},\infty)}(L_\varepsilon(T;x)) \times \int_0^{T\wedge(1/\tilde{q})} [1 - \mathbf{1}_{(-\infty,\hat{q}y/\tilde{q}]}(L_\varepsilon(s;x))]\mathbf{1}_{(x,\infty)}(\xi_\varepsilon(1-\tilde{q}s))\, ds\Big\}$$

$$= \mathbf{E}\Big\{\int_0^{T\wedge(1/\tilde{q})} \mathbf{1}_{(\hat{q}y/\tilde{q},\infty)}(L_\varepsilon(s;x))\mathbf{1}_{(x,\infty)}(\xi(1-\tilde{q}s))\, ds\Big\}$$

$$= \int_0^{T\wedge(1/\tilde{q})} \mathbf{P}\{\tilde{q}L_\varepsilon(s;x)/\hat{q} > y, \xi_\varepsilon(1-\tilde{q}s) > x\}\, ds.$$

Using (2.1) and (2.2), we get density convergence

$$\lim_{\varepsilon\downarrow 0} \frac{\mathbf{P}\{\xi_\varepsilon(1-\tilde{q}s) > x\}}{\mathbf{E}\{L_\varepsilon(T;x)\}} = \frac{\overline{G}(s)}{E(T;\overline{G})} \qquad \text{for } s \in [0,T).$$

Hence (4.2) together with Scheffé's theorem [Scheffé (1947)] show that, as $\varepsilon \downarrow 0$,

$$\int_y^\infty \frac{\mathbf{P}\{\tilde{q}L_\varepsilon(T;x) > \hat{q}z\}}{(\tilde{q}/\hat{q})\mathbf{E}\{L_\varepsilon(T;x)\}}\, dz$$

$$= \int_0^{T\wedge(1/\tilde{q})} \mathbf{P}\Big\{\int_0^{s\tilde{q}/\hat{q}} \mathbf{1}_{(x,\infty)}(\xi_\varepsilon(1-\tilde{q}s+\hat{q}t))\, dt > y \Big| \xi_\varepsilon(1-\tilde{q}s) > x\Big\} \times \frac{\mathbf{P}\{\xi_\varepsilon(1-\tilde{q}s) > x\}}{\mathbf{E}\{L_\varepsilon(T;x)\}}\, ds \tag{4.3}$$

$$\sim \int_0^{T\wedge(1/\tilde{q})} \mathbf{P}\Big\{\int_0^{s\tilde{q}/\hat{q}} \mathbf{1}_{(x,\infty)}(\xi_\varepsilon(1-\tilde{q}s+\hat{q}t))\, dt > y \Big| \xi_\varepsilon(1-\tilde{q}s) > x\Big\} \times \frac{\overline{G}(s)}{E(T;\overline{G})}\, ds.$$

Picking a $\widehat{T} > 0$, (2.3) gives convergence of moments

$$\mathbf{E}\Big\{\Big(\int_0^{(s\tilde{q}/\hat{q})\wedge\widehat{T}} \mathbf{1}_{(x,\infty)}(\xi_\varepsilon(1-\tilde{q}s+\hat{q}t))\, dt\Big)^n \Big| \xi_\varepsilon(1-\tilde{q}s) > x\Big\}$$

$$= \int_{[0,(s\tilde{q}/\hat{q})\wedge\widehat{T}]^n} \mathbf{P}\left\{\bigcap_{i=1}^{n}\{\xi_\varepsilon(1-\tilde{q}s+\hat{q}t_i)>x\}\Big|\xi_\varepsilon(1-\tilde{q}s)>x\right\}dt$$

$$\to \int_{[0,(sQ_2)\wedge\widehat{T}]^n} P_s(t)\,dt \qquad \text{as } \varepsilon\downarrow 0 \text{ for } n\in\mathbb{N}.$$

Here we have, with obvious notation,

$$\left(\int_0^{(s\tilde{q}/\hat{q})\wedge\widehat{T}} \mathbf{1}_{(x,\infty)}(\xi_\varepsilon(1-\tilde{q}s+\hat{q}t))\,dt\Big|\xi_\varepsilon(1-\tilde{q}s)>x\right)$$

$$\xrightarrow{d} \text{some random variable } \zeta_{s,\widehat{T}}$$

(since the left-hand side is bounded by $\widehat{T}$). Using this in (4.3), we get

$$\begin{aligned}
&\liminf_{\varepsilon\downarrow 0}\int_y^\infty \frac{\mathbf{P}\{\tilde{q}L_\varepsilon(T;x)>\hat{q}z\}}{(\tilde{q}/\hat{q})\mathbf{E}\{L_\varepsilon(T;x)\}}\,dz\\
&\quad\ge \int_0^{T\wedge(1/\tilde{q})} \mathbf{P}\{\zeta_{s,\widehat{T}}>y\}\frac{\overline{G}(s)}{E(T;\overline{G})}\,ds \qquad (4.4)\\
&\quad\equiv \Lambda_{\widehat{T}}(T;y)\uparrow\Lambda(T;y) \qquad \text{as } \widehat{T}\uparrow\infty.
\end{aligned}$$

On the other hand, using (2.4) if $Q_2=\infty$, by (4.3) and Markov's inequality,

$$\begin{aligned}
&\limsup_{\varepsilon\downarrow 0}\int_y^\infty \frac{\mathbf{P}\{\tilde{q}L_\varepsilon(T;x)>\hat{q}z\}}{(\tilde{q}/\hat{q})\mathbf{E}\{L_\varepsilon(T;x)\}}\,dz\\
&\quad\le \limsup_{\varepsilon\downarrow 0}\int_0^{T\wedge(1/\tilde{q})} \mathbf{P}\left\{\int_{(s\tilde{q}/\hat{q})\wedge\widehat{T}}^{s\tilde{q}/\hat{q}} \mathbf{1}_{(x,\infty)}(\xi_\varepsilon(1-\tilde{q}s+\hat{q}t))\,dt>\delta \right.\\
&\qquad\qquad\qquad\left.\Big|\xi_\varepsilon(1-\tilde{q}s)>x\right\}\\
&\qquad\qquad \times\frac{\overline{G}(s)}{E(T;\overline{G})}\,ds\\
&\qquad + \limsup_{\varepsilon\downarrow 0}\int_0^{T\wedge(1/\tilde{q})} \mathbf{P}\{\zeta_{s,\widehat{T}}>y-\delta\}\frac{\overline{G}(s)}{E(T;\overline{G})}\,ds \qquad (4.5)\\
&\quad\le \limsup_{\varepsilon\downarrow 0}\int_0^{T\wedge(1/\tilde{q})}\int_{(s\tilde{q}/\hat{q})\wedge\widehat{T}}^{s\tilde{q}/\hat{q}} \frac{\mathbf{P}\{\xi_\varepsilon(1-\tilde{q}s+\hat{q}t)>x|\xi_\varepsilon(1-\tilde{q}s)>x\}}{\delta}\,dt\\
&\qquad\qquad \times\frac{\overline{G}(s)}{E(T;\overline{G})}\,ds\\
&\qquad + \Lambda_{\widehat{T}}(T;y-\delta)\\
&\quad\to 0+\Lambda(T;y-0) \qquad \text{as } \widehat{T}\to\infty \text{ and } \delta\downarrow 0 \text{ (in that order).}
\end{aligned}$$

By (4.5) and elementary arguments, we have, for $T \in I$ and $x \in J$,

$$\begin{aligned}
&\liminf_{\varepsilon\downarrow 0} \frac{\mathbf{P}\{\sup_{t\in[0,T\wedge(1/\tilde q)]}\xi_\varepsilon(1-t\tilde q)>x\}}{(\tilde q/\hat q)\mathbf{E}\{L_\varepsilon(T;x)\}}\\
&\qquad \geq \frac{1}{y}\left[1-\limsup_{\varepsilon\downarrow 0}\int_y^\infty \frac{\mathbf{P}\{\tilde q L_\varepsilon(T;x)>\hat q z\}}{(\tilde q/\hat q)\mathbf{E}\{L_\varepsilon(T;x)\}}\,dz\right] \qquad (4.6)\\
&\qquad \to \limsup_{y\downarrow 0}\frac{1-\Lambda(T;y)}{y} \qquad \text{as } y\downarrow 0
\end{aligned}$$

in a suitable way. Here the left limit is strictly positive, since by a version of (4.5),

$$\begin{aligned}
&\limsup_{\varepsilon\downarrow 0}\int_{y+x}^\infty \frac{\mathbf{P}\{\tilde q L_\varepsilon(T;x)>\hat q z\}}{(\tilde q/\hat q)\mathbf{E}\{L_\varepsilon(T;x)\}}\,dz\\
&\quad\leq \limsup_{\varepsilon\downarrow 0}\int_0^{T\wedge(1/\tilde q)}\mathbf{P}\left\{\int_0^{(s\tilde q/\hat q)\wedge\widehat T}\mathbf{1}_{(x,\infty)}(\xi_\varepsilon(1-\tilde q s+\hat q t))\,dt>y\ \middle|\ \xi_\varepsilon(1-\tilde q s)>x\right\}\\
&\qquad\qquad \times\frac{\overline G(s)}{E(T;\overline G)}\,ds\\
&\qquad + \limsup_{\varepsilon\downarrow 0}\int_0^{T\wedge(1/\tilde q)}\int_{(s\tilde q/\hat q)\wedge\widehat T}^{s\tilde q/\hat q}\mathbf{P}\{\xi_\varepsilon(1-\tilde q s+\hat q t)>x|\xi_\varepsilon(1-\tilde q s)>x\}\,dt\\
&\qquad\qquad \times\frac{\overline G(s)}{E(T;\overline G)}\,ds\\
&\quad\leq 0+\frac{1}{2} \qquad \text{for } \widehat T\in(0,y] \text{ large enough.}
\end{aligned}$$

On the other hand, we have, for $\tilde x<\hat x<x$ in $J$ and $y>0$, by (4.4), as $\varepsilon\downarrow 0$,

$$\begin{aligned}
&\frac{\mathbf{P}\{\sup_{t\in[0,T\wedge(1/\tilde q)]}\xi_\varepsilon(1-t\tilde q)>x\}}{(\tilde q/\hat q)\mathbf{E}\{L_\varepsilon(T;\tilde x)\}}\\
&\quad\leq\frac{1}{y(\tilde q/\hat q)\mathbf{E}\{L_\varepsilon(T;\tilde x)\}}\\
&\qquad\times\int_0^y\mathbf{P}\left\{\left\{\sup_{t\in[0,T\wedge(1/\tilde q)]}\xi_\varepsilon(1-t\tilde q)>x\right\}\cup\{\tilde q L_\varepsilon(T;\tilde x)>\hat q z\}\right.\\
&\qquad\qquad\left.\cup\bigcup_{k=0}^{\lfloor((T\tilde q)\wedge 1)/(a\hat q)\rfloor}\{\xi_\varepsilon(1-ka\hat q)>\hat x\}\right\}dz
\end{aligned}$$

$$
\begin{aligned}
&\le \frac{1}{y}\left[1 - \int_y^\infty \frac{\mathbf{P}\{\tilde q L_\varepsilon(T;\tilde x) > \hat q z\}}{(\tilde q/\hat q)\mathbf{E}\{L_\varepsilon(T;\tilde x)\}}\,dz\right] \\
&\quad + \sum_{k=0}^{\lfloor((T\tilde q)\wedge 1)/(a\hat q)\rfloor} \frac{\mathbf{P}\{\tilde q L_\varepsilon(T;\tilde x) \le \hat q y, \xi_\varepsilon(1-ka\hat q) > \hat x\}}{(\tilde q/\hat q)\mathbf{E}\{L_\varepsilon(T;\tilde x)\}} \\
&\quad + \frac{1}{(\tilde q/\hat q)\mathbf{E}\{L_\varepsilon(T;\tilde x)\}}\mathbf{P}\Bigg\{\sup_{t\in[0,T\wedge(1/\tilde q)]} \xi_\varepsilon(1-t\tilde q) > x, \\
&\qquad\qquad \bigcap_{k=0}^{\lfloor((T\tilde q)\wedge 1)/(a\hat q)\rfloor} \{\xi_\varepsilon(1-ka\hat q) \le \hat x\}\Bigg\} \qquad (4.7) \\
&\le \frac{1-\Lambda(T;y-0)}{y} + o(1) \\
&\quad + \frac{\mathbf{P}\{\tilde q L_\varepsilon(T;\tilde x) \le \hat q y, \xi_\varepsilon(1) > \hat x\}}{(\tilde q/\hat q)\mathbf{E}\{L_\varepsilon(T;\tilde x)\}} \\
&\quad + \sum_{k=1}^{\lfloor((T\tilde q)\wedge 1)/(a\hat q)\rfloor} \int_{k-1}^{k} [\mathbf{P}\{\tilde q L_\varepsilon(T;\tilde x) \le \hat q y, \\
&\qquad\qquad \xi_\varepsilon(1-ta\hat q) \vee \xi_\varepsilon((1-(t+1)a\hat q)^+) > \tilde x\} \\
&\qquad\qquad \times ((\tilde q/\hat q)\mathbf{E}\{L_\varepsilon(T;\tilde x)\})^{-1} \\
&\qquad\qquad + \mathbf{P}\{\xi_\varepsilon(1-ka\hat q) > \hat x, \\
&\qquad\qquad\quad \xi_\varepsilon(1-ta\hat q) \vee \xi_\varepsilon((1-(t+1)a\hat q)^+) \le \tilde x\} \\
&\qquad\qquad\qquad \times ((\tilde q/\hat q)\mathbf{E}\{L_\varepsilon(T;\tilde x)\})^{-1}]\,dt \\
&\quad + \frac{1}{(\tilde q/\hat q)\mathbf{E}\{L_\varepsilon(T;\tilde x)\}} \\
&\quad\times \sum_{k=0}^{\lfloor((T\tilde q)\wedge 1)/(a\hat q)\rfloor} \mathbf{P}\Bigg\{\sup_{t\in[k,k+1]} \xi_\varepsilon(1-ta\tilde q) > x, \\
&\qquad\qquad \xi_\varepsilon(1-ka\hat q) \vee \xi_\varepsilon((1-(k+1)a\hat q)^+) \le \hat x\Bigg\} \\
&\le \frac{1-\Lambda(T;y-0)}{y} + o(1)
\end{aligned}
$$

$$
\begin{aligned}
&+\limsup_{\varepsilon\downarrow 0}\frac{\mathbf{P}\{\tilde q L_\varepsilon(T;\tilde x)\le \hat q y,\xi_\varepsilon(1)>\hat x\}}{(\tilde q/\hat q)\mathbf{E}\{L_\varepsilon(T;\tilde x)\}}\\
&+\limsup_{\varepsilon\downarrow 0} 2\int_0^{\lfloor T\wedge(1/\tilde q)\rfloor}\frac{\mathbf{P}\{\tilde q L_\varepsilon(t;\tilde x)\le \hat q y,\xi_\varepsilon(1-t\tilde q)>\tilde x\}}{a\mathbf{E}\{L_\varepsilon(T;\tilde x)\}}\,dt\\
&+\limsup_{\varepsilon\downarrow 0}\frac{1}{a\mathbf{E}\{L_\varepsilon(T;\tilde x)\}}\\
&\qquad\times\int_0^{T\wedge(1/\tilde q)}\mathbf{P}\bigg\{\sup_{s\in[0,1]}\xi_\varepsilon((1-t\tilde q-sa\hat q)^+)>\hat x,\\
&\qquad\qquad\qquad \xi_\varepsilon(1-t\tilde q)\vee\xi_\varepsilon((1-t\tilde q-a\hat q)^+)\le\tilde x\bigg\}\,dt\\
&+\limsup_{\varepsilon\downarrow 0}\frac{1}{a\mathbf{E}\{L_\varepsilon(T;\tilde x)\}}\\
&\qquad\times\int_0^{T\wedge(1/\tilde q)}\mathbf{P}\bigg\{\sup_{s\in[0,1]}\xi_\varepsilon((1-t\tilde q-sa\hat q)^+)>x,\\
&\qquad\qquad\qquad \xi_\varepsilon(1-t\tilde q)\vee\xi_\varepsilon((1-t\tilde q-a\hat q)^+)\le\hat x\bigg\}\,dt.
\end{aligned}
$$

Here lim sup of the left-hand side is finite, since the first term on the right-hand side is trivially finite, the second by (4.1) and the weak convergence used in (4.4), the third by (2.2) and (4.1) and the fourth and fifth, for $a>0$ small enough, by (2.1), (4.1) and (2.5). Hence the right-hand side limit in (4.6) is finite. Sending $y\downarrow 0$ and $\tilde x\uparrow x$ in (4.7), and using (2.1), (4.1) and (2.5)–(2.6), we see that the following two limits exist and coincide:

$$
\lim_{\varepsilon\downarrow 0}\frac{\mathbf{P}\{\sup_{t\in[0,T\wedge(1/\tilde q)]}\xi_\varepsilon(1-t\tilde q)>x\}}{(\tilde q/\hat q)\mathbf{E}\{L_\varepsilon(T;x)\}}=\lim_{y\downarrow 0}\frac{1-\Lambda(T;y)}{y}\equiv\kappa(T)>0. \tag{4.8}
$$

In addition, it follows that $\kappa\equiv\sup_{T\in I}\kappa(T)<\infty$. Moreover, by inspection of (4.7), for $T\in I$ and $\hat x<x$ in $J$,

$$
\begin{aligned}
&\lim_{a\downarrow 0}\limsup_{\varepsilon\downarrow 0}\frac{1}{q}\mathbf{P}\bigg\{\sup_{t\in[0,T\wedge(1/\tilde q)]}\xi_\varepsilon(1-t\tilde q)>x,\\
&\qquad\qquad \bigcap_{k=0}^{\lfloor((T\tilde q)\wedge 1)/(a\hat q)\rfloor}\{\xi_\varepsilon(1-ka\hat q)\le\hat x\}\bigg\}=0.
\end{aligned}\tag{4.9}
$$

By (2.7), (4.8) and (4.9), we have, for $\hat t\in[0,1)$, $\hat x<x$ in $J$ and $T\in I$,

$$
\liminf_{h\downarrow 0}\liminf_{\varepsilon\downarrow 0}\frac{1}{h}\mathbf{P}\bigg\{\bigcup_{n=\lfloor\hat t/q\rfloor+1}^{\lfloor(\hat t+h)/q\rfloor}\ \bigcup_{k=0}^{\lfloor((T\tilde q)\wedge 1)/(a\hat q)\rfloor}\{\xi_\varepsilon(n-ka\hat q)>\hat x\}\bigg\}
$$

$$\geq \liminf_{h\downarrow 0}\liminf_{\varepsilon\downarrow 0}\frac{\lfloor h/q\rfloor}{h}\mathbf{P}\Bigg\{\bigcup_{k=0}^{\lfloor((T\tilde{q})\wedge 1)/(a\hat{q})\rfloor}\{\xi_\varepsilon(1-ka\hat{q})>\hat{x}\}\Bigg\}$$

$$-\limsup_{h\downarrow 0}\limsup_{\varepsilon\downarrow 0}\frac{1}{h}\sum_{m=\lfloor\hat{t}/q\rfloor+1}^{\lfloor(\hat{t}+h)/q\rfloor}\sum_{n=m+1}^{\lfloor(\hat{t}+h)/q\rfloor}\sum_{k,\ell=0}^{\lfloor((T\tilde{q})\wedge 1)/(a\hat{q})\rfloor}$$

$$\mathbf{P}\{\xi_\varepsilon(m-ka\hat{q})>\hat{x},\xi_\varepsilon(n-\ell a\hat{q})>\hat{x}\}$$

$$\geq \kappa(T)\overline{F}(x)+o(1) \tag{4.10}$$

$$-\limsup_{h\downarrow 0}\limsup_{\varepsilon\downarrow 0}$$

$$\sup_{r,s\in[0,T\wedge(1/\tilde{q}))}\sum_{n=2}^{\lfloor h/q\rfloor}\frac{q\mathbf{P}\{\xi_\varepsilon(1-r\tilde{q})>\hat{x},\xi_\varepsilon(n-s\tilde{q})>\hat{x}\}}{\mathbf{P}\{\xi_\varepsilon(1-r\tilde{q})>0\}\mathbf{P}\{\xi_\varepsilon(1-s\tilde{q})>0\}}$$

$$\times\Bigg(\sum_{k=0}^{\lfloor((T\tilde{q})\wedge 1)/(a\hat{q})\rfloor}\frac{\mathbf{P}\{\xi_\varepsilon(1-ka\hat{q})>0\}}{(\hat{q}/\tilde{q})\mathbf{P}\{\xi_\varepsilon(1)>0\}}\Bigg)^2$$

$$=\kappa(T)\overline{F}(x)+o(1)\to\kappa(T)\overline{F}(x)\qquad\text{as } a\downarrow 0,$$

since (2.1), (2.2) and (2.5) show that the sum on the last row is bounded as $\varepsilon\downarrow 0$.

To finish the proof, notice that, given $T\in I$ and $x\in J$, (2.8) yields

$$\begin{aligned}&\lim_{\varepsilon\downarrow 0}\Bigg[\mathbf{P}\Bigg\{\bigcap_{m=1}^{\lfloor 1/h\rfloor}\bigcap_{n=\lfloor(m-1)h/q\rfloor+1}^{\lfloor(mh-h^2)/q\rfloor}\bigcap_{k=0}^{\lfloor((T\tilde{q})\wedge 1)/(a\hat{q})\rfloor-1}\{\xi_\varepsilon(n-ka\hat{q})\leq x\}\Bigg\}\\&\qquad-\mathbf{P}\Bigg\{\bigcap_{n=1}^{\lfloor(h-h^2)/q\rfloor}\bigcap_{k=0}^{\lfloor((T\tilde{q})\wedge 1)/(a\hat{q})\rfloor-1}\{\xi_\varepsilon(n-ka\hat{q})\leq x\}\Bigg\}^{\lfloor 1/h\rfloor}\Bigg]=0.\end{aligned} \tag{4.11}$$

Using (4.11) together with (4.10), we readily obtain, for $x\in J$,

$$\limsup_{\varepsilon\downarrow 0}\mathbf{P}\Big\{\sup_{r\in[0,1/q]}\xi_\varepsilon(r)\leq x\Big\}$$

$$\leq\inf_{T\in I}\limsup_{a\downarrow 0}$$

$$\limsup_{h\downarrow 0}\Bigg(1-h\limsup_{\varepsilon\downarrow 0}\frac{1}{h}\mathbf{P}\Bigg\{\bigcup_{n=1}^{\lfloor(h-h^2)/q\rfloor}\bigcup_{k=0}^{\lfloor((T\tilde{q})\wedge 1)/(a\hat{q})\rfloor-1} \tag{4.12}$$

$$\left.\left.\{\xi_\varepsilon(n-ka\hat{q})>x\}\right\}\right)^{\lfloor 1/h\rfloor}$$

$$\leq \inf_{T\in I} e^{-\kappa(T)\overline{F}(\tilde{x})} = e^{-\kappa\overline{F}(\tilde{x})} \to e^{-\kappa\overline{F}(x)} \qquad \text{as } \tilde{x}\downarrow x.$$

On the other hand, (4.11) together with (4.8) and (4.9) similarly gives, for $x\in J$,

$$\liminf_{\varepsilon\downarrow 0}\mathbf{P}\left\{\sup_{r\in[0,1/q]}\xi_\varepsilon(r)\leq x\right\}$$

$$\geq \inf_{T\in I}\liminf_{a\downarrow 0}$$

$$\text{(4.13)}\qquad \liminf_{h\downarrow 0}\liminf_{\varepsilon\downarrow 0}\mathbf{P}\left\{\bigcap_{m=1}^{\lfloor 1/h\rfloor}\bigcap_{n=\lfloor(m-1)h/q\rfloor+1}^{\lfloor(mh-h^2)/q\rfloor}\bigcap_{k=0}^{\lfloor((T\tilde{q})\wedge 1)/(a\hat{q})\rfloor-1}\{\xi_\varepsilon(n-ka\hat{q})\leq\hat{x}\}\right\}$$

$$\geq \inf_{T\in I}\liminf_{a\downarrow 0}$$

$$\liminf_{h\downarrow 0}\liminf_{\varepsilon\downarrow 0}\left(1-\frac{h}{q}\mathbf{P}\left\{\bigcup_{k=0}^{\lfloor((T\tilde{q})\wedge 1)/(a\hat{q})\rfloor}\{\xi_\varepsilon(n-ka\tilde{q})>\hat{x}\}\right\}\right)^{\lfloor 1/h\rfloor}$$

$$\geq \inf_{T\in I} e^{-\kappa(T)\overline{F}(\hat{x})} = e^{-\kappa\overline{F}(\hat{x})} \to e^{-\kappa\overline{F}(x)} \qquad \text{as } \hat{x}\uparrow x.$$

Clearly, (4.12) and (4.13) add up to (4.1) or (4.2), depending on whether $i=1$ or 2. □

**5. Stable processes.** We give two examples with strictly *α-stable Lévy motions*, and prove (1.1) for totally skewed *linear fractional α-stable motion* (LFSM).

Let $\{L(t)\}_{t\in\mathbb{R}}$ be a separable $\alpha$-stable Lévy motion with $L(t)\sim S_\alpha(|t|^{1/\alpha},\beta,\mu)$, where $\alpha\in(0,2)$, $\beta\in[-1,1]$ and $\mu\in\mathbb{R}$, with $\mu=0$ if $\alpha\neq 1$ and $\beta=0$ if $\alpha=1$. [See Samorodnitsky and Taqqu (1994) on basic facts and notation for $\alpha$-stable processes.]

It is tempting to try to extend Example 3 to a general SS Lévy process $L$:

EXAMPLE 4. For $\beta\neq -1$, Willekens (1987) showed that [see also Berman (1986)]

$$\text{(5.1)}\ \mathbf{P}\left\{\sup_{t\in[0,1]}L(t)>u\right\}\sim\mathbf{P}\{L(1)>u\}\sim\tfrac{1}{2}C_\alpha(1+\beta)u^{-\alpha} \qquad \text{as } u\to\infty,$$

where $C_\alpha = (\int_0^\infty x^{-\alpha} \sin(x)\, dx)^{-1}$. Hence, since $L$ is SS with index $1/\alpha(1/\alpha$-SS),

$$\begin{aligned}
&\mathbf{P}\Big\{\sup_{t\in[0,1]} [L(t) - L_q(t)] \le \varepsilon\Big\} \\
&\quad = \Big(\mathbf{P}\Big\{\sup_{t\in[0,q]} L(t) \le \varepsilon\Big\}\Big)^{\lfloor 1/q\rfloor} \mathbf{P}\Big\{\sup_{t\in[0,1-\lfloor 1/q\rfloor q]} L(t) \le \varepsilon\Big\} \\
&\quad \sim (1 - \mathbf{P}\{L(1) > \varepsilon/q^{1/\alpha}\})^{\lfloor 1/q\rfloor} \\
&\quad \to \exp\{-\tfrac{1}{2} C_\alpha (1+\beta)\varepsilon^{-\alpha}\} \qquad \text{as } q \downarrow 0 \text{ for } \varepsilon > 0!
\end{aligned}$$

[It is an exercise in subexponentiality to use this in turn to recover (5.1).] Here $L$ has jumps with positive probability. Hence the sampling problem does not really make sense: we are in fact computing the distribution of the largest positive jump.

Further, again by (5.1), and in the case when $\alpha > 1$ and $\beta = -1$ also using that

$$\mathbf{P}\Big\{\sup_{t\in[0,h]} L(t) > u\Big\} \le \frac{\mathbf{P}\{L(h) > u\}}{\mathbf{P}\{L(h) > 0\}} \qquad \text{for } h, u > 0 \tag{5.2}$$

[e.g., Albin (1993), equation 1.2], together with (5.3), we similarly get

$$\begin{aligned}
&(1 - \mathbf{P}\{|L(1)| > \varepsilon/q^{1/\alpha}\})^{\lfloor 1/q\rfloor} \\
&\quad \ge \mathbf{P}\Big\{\sup_{t\in[0,1]} |L(t) - L_q(t)| \le \varepsilon\Big\} \\
&\quad \ge \Big(1 - \mathbf{P}\Big\{\sup_{t\in[0,1]} L(t) > \varepsilon/q^{1/\alpha}\Big\} - \mathbf{P}\Big\{\inf_{t\in[0,1]} L(t) < -\varepsilon/q^{1/\alpha}\Big\}\Big)^{1+1/q} \\
&\quad \sim (1 - \mathbf{P}\{|L(1)| > \varepsilon/q^{1/\alpha}\})^{\lfloor 1/q\rfloor} \\
&\quad \to \exp\{-C_\alpha \varepsilon^{-\alpha}\} \qquad \text{as } q \downarrow 0.
\end{aligned}$$

For one case not covered by Example 4, that is, (1.1) for $\beta = -1$, Theorem 1 applies:

EXAMPLE 5. Let $\alpha > 1$ and $\beta = -1$. Denote

$$q(\varepsilon) \equiv \frac{\varepsilon^\alpha}{b_\alpha^{\alpha/\lambda_\alpha}} \Big[2\alpha \ln\Big(\frac{1}{\varepsilon}\Big) - (3 - 2\alpha)\ln\Big(2\alpha\ln\Big(\frac{1}{\varepsilon}\Big)\Big) - 2\ln\Big(\frac{\sqrt{2\pi\alpha}}{b_\alpha^{\alpha/\lambda_\alpha}}\Big)\Big]^{-\alpha/(2\lambda_\alpha)},$$

where

$$b_\alpha^2 = \frac{\alpha^{2\lambda_\alpha}}{2(\alpha - 1)|\cos(\pi\alpha/2)|^{2\lambda_\alpha - 1}} \quad \text{and} \quad \lambda_\alpha = \frac{\alpha}{2(\alpha - 1)}.$$

Recall that, as $u \to \infty$ [e.g., Samorodnitsky and Taqqu (1994), equation 1.2.11],

$$\begin{aligned}\mathbf{P}\{S_\alpha(\sigma,-1,0)>u\} &= \mathbf{P}\{S_\alpha(1,-1,0)>u/\sigma\} \\ &\sim \frac{\mathbf{P}\{\mathrm{N}(0,b_\alpha^2)>(u/\sigma)^{\lambda_\alpha}\}}{\sqrt{\alpha}}.\end{aligned} \tag{5.3}$$

Taking $w \equiv \varepsilon/(2\alpha\lambda_\alpha \ln(1/\varepsilon))$ and $\tilde{q} = \hat{q} \equiv \alpha w/\varepsilon$ ($Q_1 = \infty$ and $Q_2 = 1$), (5.3) gives (2.1) and (2.2) for $L_\varepsilon = L_{1,\varepsilon}$, with $\overline{F} = \overline{G} = e^{-\cdot}$, $J = \mathbb{R}$ and $I = (0,\infty)$.

Let $\zeta$ be a standard exponential random variable independent of $L$. Since

$$(L_\varepsilon(1-\tilde{q}r)|L_\varepsilon(1-\tilde{q}r)>x) \xrightarrow{d} \zeta + x,$$

by (2.1), we now get (2.3) by observing that, as $\varepsilon \downarrow 0$,

$$\begin{aligned}&\{(L_\varepsilon(1+\tilde{q}(t-r))|L_\varepsilon(1-\tilde{q}r)>x)\}_{t\in(0,r]} \\ &\quad \stackrel{d}{=} \left\{\frac{L(q\tilde{q}t)}{w} + (L_\varepsilon(1-\tilde{q}r)|L_\varepsilon(1-\tilde{q}r)>x)\right\}_{t\in(0,r]} \\ &\quad \xrightarrow{d} \left\{L\left(\left|\cos\left(\frac{\pi\alpha}{2}\right)\right|t\right) + x + \zeta\right\}_{t\in(0,r]}.\end{aligned}$$

Using (2.1) together with (5.2), we get (2.6) in the following way:

$$\begin{aligned}&\limsup_{\varepsilon\downarrow 0} \frac{\mathbf{P}\{\sup_{r\in[T\wedge(1/\tilde{q}),(1/\tilde{q}))} L_\varepsilon(1-r\tilde{q})>x\}}{\mathbf{P}\{L_\varepsilon(1)>0\}} \\ &\quad \le \limsup_{\varepsilon\downarrow 0} \frac{\mathbf{P}\{L_\varepsilon(1-T\tilde{q})>x\}}{\mathbf{P}\{L(q-T\tilde{q}q)>0\}\mathbf{P}\{L_\varepsilon(1)>0\}} \\ &\quad = \frac{e^{-T-x}}{\mathbf{P}\{L(1)>0\}} \to 0 \qquad \text{as } T\to\infty.\end{aligned}$$

Using $1/\alpha$-SS and (5.2), (5.3), (2.5) follows from the fact that, for $y < y+\delta$,

$$\begin{aligned}&\frac{1}{a\mathbf{P}\{L_\varepsilon(1)>0\}}\mathbf{P}\left\{\sup_{s\in[0,a]} L_\varepsilon((1-r\tilde{q}-s\tilde{q})^+)>y+\delta, L_\varepsilon((1-r\tilde{q}-a\tilde{q})^+)\le y\right\} \\ &\quad \le \frac{1}{a\mathbf{P}\{L_\varepsilon(1)>0\}}\sum_{k=0}^{\infty}\mathbf{P}\left\{\sup_{s\in[0,a]} L_\varepsilon((1-r\tilde{q}-s\tilde{q})^+)\right. \\ &\qquad\qquad - L_\varepsilon((1-r\tilde{q}-a\tilde{q})^+)>\delta+k,\end{aligned}$$

$$L_\varepsilon((1-r\tilde{q}-a\tilde{q})^+) > y-k-1\Big\}$$

$$\leq \sum_{k=0}^{\infty} \frac{\mathbf{P}\{L(a\tilde{q}q) > (\delta+k)w\}\mathbf{P}\{L_\varepsilon(1) > y-k-1\}}{a\mathbf{P}\{L(1)>0\}\mathbf{P}\{L_\varepsilon(1)>0\}}$$

$$\leq \sum_{k=0}^{\infty} \frac{\mathbf{P}\{L(|\cos((\pi\alpha)/2)|) > (1/2)(\delta+k)/a^{1/\alpha}\}\mathbf{P}\{L_\varepsilon(1) > y-k-1\}}{a\mathbf{P}\{L(1)>0\}\mathbf{P}\{L_\varepsilon(1)>0\}}$$

$$\to 0 \qquad \text{uniformly in } r$$

as $\varepsilon \downarrow 0$ and $a \downarrow 0$, since, by (5.3),

$$\frac{\mathbf{P}\{L_\varepsilon(1) > y-k-1\}}{\mathbf{P}\{L(1)>0\}} \leq K_y e^{2k}.$$

Finally, (2.7) and (2.8) follow trivially from independence of increments.

Now take $\alpha > 1$ and $\beta = -1$, as in Example 5. Given an $H \in (1/\alpha, 1)$, consider the LFSM given by

$$\{\xi(t)\}_{t\geq 0} \stackrel{d}{=} \left\{\int_{r\in\mathbb{R}} [((t+r)^+)^{H-1/\alpha} - (r^+)^{H-1/\alpha}]\,dL(r)\right\}_{t\geq 0}. \tag{5.4}$$

This process is $H$-SSSI and totally skewed to the left, $\xi(t) \sim S_\alpha(\sigma(\xi(1))t^H, -1, 0)$.

THEOREM 2. *Let $\xi$ be given by* (5.4). *With*

$$q(\varepsilon) \equiv \left(\frac{\varepsilon}{b_\alpha^{1/\lambda_\alpha}\sigma(\xi(1))}\right)^{1/H} \times \left[\frac{2}{H}\ln\left(\frac{1}{\varepsilon}\right) - \frac{H\lambda_\alpha - 1}{H\lambda_\alpha}\ln\left(\frac{2}{H}\ln\left(\frac{1}{\varepsilon}\right)\right) - 2\ln\left(\frac{\sqrt{2\pi\alpha}}{[b_\alpha^{1/\lambda_\alpha}\sigma(\xi(1))]^{1/H}}\right)\right]^{-1/(2H\lambda_\alpha)}$$

*and $w(\varepsilon) \equiv H\varepsilon/(2\lambda_\alpha \ln(1/\varepsilon))$ (cf. Example 4), we have*

$$\lim_{\varepsilon\downarrow 0} \mathbf{P}\left\{\frac{\sup_{t\in[0,1]}[\xi(t) - \xi(\lfloor t/q(\varepsilon)\rfloor q(\varepsilon))] - \varepsilon}{w(\varepsilon)} \leq x\right\} = \exp\{-e^{-x}\} \qquad \textit{for } x \in \mathbb{R}.$$

PROOF. With $\tilde{q}(\varepsilon) = \hat{q}(\varepsilon) = w(\varepsilon)/(H\varepsilon)$, so that $Q_1 = \infty$ and $Q_2 = 1$, (5.3) gives (2.1) and (2.2), for $\overline{F} = \overline{G} = e^{-\cdot}$, $J = \mathbb{R}$ and $I = (0, \infty)$, as in Example 5.

By Albin [(1998), Section 13], we have, with $v(u) \sim c_1 u^{-1/(\alpha-1)}$ and $p(u) \sim c_2 v(u)^{1/H}$ as $u \to \infty$, for some constants $c_1, c_2 > 0$,

$$\begin{aligned} &\mathbf{P}\{\xi(1-p(u)t) > u + (\delta + x)v(u), \xi(1) \le u + xv(u)\} \\ &\qquad \le C_1 \delta^{-d} t^e \mathbf{P}\{\xi(1) > u + xv(u)\} \end{aligned} \tag{5.5}$$

for $u$ large enough, $t > 0$ small enough and $\delta > 0$, for some constants $C_1, d > 0$ and $e > 1$. For any function $\hat{p}(u) = o(p(u))$ as $u \to \infty$, we thus have, by (2.1),

$$\lim_{u \to \infty} \frac{\mathbf{P}\{\xi(1-\hat{p}(u)t) > u + (\delta + x)v(u), \xi(1) \le u + xv(u)\}}{\mathbf{P}\{\xi(1) > u\}} = 0.$$

From this we get, by (2.1) and $H$-SS, with the change of variable $u = \varepsilon/q^H$, $v(u) = w/q^H \sim c_1 u^{-1/(\alpha-1)}$ and $\hat{p}(u) = \tilde{q}/(1 - \tilde{q}(r-t)) = o(v(u)^{1/H})$, for $t \in (0, r)$,

$$\begin{aligned} &\frac{\mathbf{P}\{\xi_\varepsilon(1-\tilde{q}r) > \delta + x, \xi_\varepsilon(1-\tilde{q}(r-t)) \le x\}}{\mathbf{P}\{\xi_\varepsilon(1-\tilde{q}r) > x\}} \\ &\quad \sim \frac{e^{r+x}}{\mathbf{P}\{\xi(1) > \varepsilon/q^H\}} \\ &\qquad \times \mathbf{P}\biggl\{\xi\biggl(1 - \frac{\tilde{q}t}{1-\tilde{q}(r-t)}\biggr) > \frac{\varepsilon + w(\delta + x)}{q^H[1-\tilde{q}(r-t)]^H}, \\ &\qquad\qquad\qquad \xi(1) \le \frac{\varepsilon + wx}{q^H[1-\tilde{q}(r-t)]^H}\biggr\} \\ &\quad \sim \frac{e^{r+x}}{\mathbf{P}\{\xi(1) > \varepsilon/q^H\}} \\ &\qquad \times \mathbf{P}\biggl\{\xi\biggl(1 - \frac{\tilde{q}t}{1-\tilde{q}(r-t)}\biggr) > \frac{\varepsilon + w(\delta + x + r - t)}{q^H}, \\ &\qquad\qquad\qquad \xi(1) \le \frac{\varepsilon + w(x + r - t)}{q^H}\biggr\} \\ &\quad = e^{r+x} \frac{\mathbf{P}\{\xi(1-\hat{p}t) > u + (\delta + x + r - t)v, \xi(1) \le u + (x + r - t)v\}}{\mathbf{P}\{\xi(1) > u\}} \\ &\quad \to 0 \qquad \text{as } u \to \infty. \end{aligned} \tag{5.6}$$

This in turn gives (2.3) with $P_r(t) = 1$ since, by (2.1),

$$\begin{aligned} &\liminf_{\varepsilon \downarrow 0} \mathbf{P}\{\xi_\varepsilon(1-\tilde{q}(r-t)) > x | \xi_\varepsilon(1-\tilde{q}r) > x\} \\ &\qquad \ge \liminf_{\delta \downarrow 0} \liminf_{\varepsilon \downarrow 0} \frac{\mathbf{P}\{\xi_\varepsilon(1-\tilde{q}(r-t)) > x, \xi_\varepsilon(1-\tilde{q}r) > \delta + x\}}{\mathbf{P}\{\xi_\varepsilon(1-\tilde{q}r) > x\}} \end{aligned}$$

$$\geq \liminf_{\delta\downarrow 0}\liminf_{\varepsilon\downarrow 0}\frac{\mathbf{P}\{\xi_\varepsilon(1-\tilde{q}r)>\delta+x\}}{\mathbf{P}\{\xi_\varepsilon(1-\tilde{q}r)>x\}}$$

$$-\limsup_{\delta\downarrow 0}\limsup_{\varepsilon\downarrow 0}\frac{\mathbf{P}\{\xi_\varepsilon(1-\tilde{q}(r-t))\leq x,\xi_\varepsilon(1-\tilde{q}r)>\delta+x\}}{\mathbf{P}\{\xi_\varepsilon(1-\tilde{q}r)>x\}}$$

$$=\liminf_{\delta\downarrow 0}e^{-\delta}-0=1.$$

By (5.5) and $H$-SS, we get easily that [this is in essence what is going on in (5.6)]

$$\mathbf{P}\{\xi(1-pt)>u+(\delta+x)v,\xi(1-ps)\leq u+xv\}\leq C_2\delta^{-d}(t-s)^e\mathbf{P}\{\xi(1)>u+xv\}$$

for $u$ large enough, $t-s>0$ small enough, $\delta>0$ and $x$ and $t$ in compacts, for some constants $C_2,d>0$ and $e>1$. Hence Albin [(1992), Proposition 2] (which does not really use stationarity) shows that, for $x<\tilde{x}$, uniformly for $s\geq 0$ in compacts,

$$\lim_{a\downarrow 0}\limsup_{u\to\infty}\frac{\mathbf{P}\{\sup_{t\in[0,a]}\xi(1-p(s+t))>u+\tilde{x}v,\xi(1-ps)\leq u+xv\}}{a\mathbf{P}\{\xi(1)>u\}}=0.$$

By the change of variable $u=\varepsilon/q^H$, and using $H$-SS, we now get (2.5), since

$$\lim_{a\downarrow 0}\limsup_{\varepsilon\downarrow 0}\frac{\mathbf{P}\{\sup_{t\in[0,a]}\xi_\varepsilon(1-\tilde{q}(s+t))>\tilde{x},\xi_\varepsilon(1-\tilde{q}s)\leq x\}}{a\mathbf{P}\{\xi_\varepsilon(1)>0\}}=0.$$

Further, (2.6) follows easily from Albin [(1998), Theorem 12], together with (2.1).

By $H$-SS, with the new variable $u$ [and since $\tilde{q}(u)\sim c_3/u$], (2.8) follows if, given $T>0$, $x\in\mathbb{R}$ and $a,h>0$, for any $s_1,\dots,s_k,t_1,\dots,t_{k'}\in[0,1/\mathbf{P}\{\xi(1)>u\}]\cap\bigcup_{n=1}^\infty(n-T/u,n]$ with $s_{i+1}-s_i\geq a/u$, $t_{j+1}-t_j\geq a/u$ and $t_1-s_k\geq h/\mathbf{P}\{\xi(1)>u\}$,

$$\lim_{u\to\infty}\Biggl[\mathbf{P}\Biggl\{\bigcap_{i=1}^{k}\{\xi(s_i)-\xi(\lfloor s_i\rfloor)\leq u+xv\},\bigcap_{j=1}^{k'}\{\xi(t_j)-\xi(\lfloor t_j\rfloor)\leq u+xv\}\Biggr\}$$

$$-\mathbf{P}\Biggl\{\bigcap_{i=1}^{k}\{\xi(s_i)-\xi(\lfloor s_i\rfloor)\leq u+xv\}\Biggr\}\times\mathbf{P}\Biggl\{\bigcap_{j=1}^{k'}\{\xi(t_j)-\xi(\lfloor t_j\rfloor)\leq u+xv\}\Biggr\}\Biggr]=0. \tag{5.7}$$

We show (5.7) by "truncating" the stochastic integral for $\xi(t)-\xi(\lfloor t\rfloor)$ from (5.4),

$$\xi(t)-\xi(\lfloor t\rfloor)=\int_{r\in\mathbb{R}}[((t+r)^+)^{H-1/\alpha}-((\lfloor t\rfloor+r)^+)^{H-1/\alpha}]\,dL(r)$$

$$= \int_{|r+t| \le (1/2)h/\mathbf{P}\{\xi(1)>u\}} + \int_{|r+t|>(1/2)h/\mathbf{P}\{\xi(1)>u\}}$$
$$\equiv I_1(t) + I_2(t).$$

Since $k + k' \le 1 + T/[a\mathbf{P}\{\xi(1) > u\}]$, and by (2.1), it is enough to show that

$$\limsup_{u\to\infty} \frac{T\mathbf{P}\{I_2 > \delta v\}}{a\mathbf{P}\{\xi(1) > u\}} = 0 \qquad \text{for } \delta > 0, \tag{5.8}$$

because then we can subtract all $I_2(s_i)$ and $I_2(t_j)$ from the corresponding process values and, after having factorized the resulting independent probabilities, add the $I_2$'s again. The approximation error is asymptotically negligible, since by (2.1),

$$\lim_{\delta\downarrow 0} \limsup_{u\to\infty} (k + k')\mathbf{P}\{\xi(1) \in [u + xv - \delta v, u + xv + \delta v]\} = 0.$$

We get (5.8) from (5.3) using that $(1 - x)^{H-1/\alpha} \ge 1 - x$ since, uniformly in $t$,

$$\begin{aligned}
\sigma(I_2(t))^\alpha &= \int_{|r+t|>(1/2)h/\mathbf{P}\{\xi(1)>u\}} [(t + r)^{H-1/\alpha} - ((\lfloor t\rfloor + r)^+)^{H-1/\alpha}]^\alpha \, dr \\
&= \int_{r+t>(1/2)h/\mathbf{P}\{\xi(1)>u\}} (t + r)^{\alpha H-1} \left(1 - \left[1 - \frac{t - \lfloor t\rfloor}{t + r}\right]^{H-1/\alpha}\right)^\alpha dr \\
&\le \int_{r+t>(1/2)h/\mathbf{P}\{\xi(1)>u\}} (t + r)^{\alpha H-1-\alpha} \, dr \\
&= O([\mathbf{P}\{\xi(1) > u\}]^{\alpha(1-H)}).
\end{aligned}$$

With the new variable $u$, using (2.1), (2.7) can be rewritten

$$\begin{aligned}
\lim_{h\downarrow 0} \limsup_{u\to\infty} \sup_{r,s\in[0,T)} \sum_{n=2}^{\lfloor h/\mathbf{P}\{\xi(1)>u\}\rfloor} & \mathbf{P}\{\xi(1 - r/u) - \xi(\lfloor 1 - r/u\rfloor) > u, \\
& \xi(n - s/u) - \xi(\lfloor n - s/u\rfloor) > u\} \\
& \times [\mathbf{P}\{\xi(1) > u\}]^{-1} = 0.
\end{aligned}$$

By positivities present, this holds if

$$\begin{aligned}
\lim_{h\downarrow 0} \limsup_{u\to\infty} \sup_{s\in[0,T)} \sum_{n=2}^{\lfloor h/\mathbf{P}\{\xi(1)>u\}\rfloor} & \mathbf{P}\{\xi(1) + \xi(n - s/u) - \xi(\lfloor n - s/u\rfloor) > 2u\} \\
& \times [\mathbf{P}\{\xi(1) > u\}]^{-1} = 0.
\end{aligned} \tag{5.9}$$

When checking (5.9), we can disregard any finite number of terms, since by Minkowski's inequality, and more or less immediate properties of LFSM, the scale

$$\begin{aligned}&\sigma(\xi(1)+\xi(n-s/u)-\xi(\lfloor n-s/u\rfloor))\\&\quad\le\sigma(\xi(1))+\sigma(\xi(n-s/u)-\xi(\lfloor n-s/u\rfloor))-\delta\\&\quad=2\sigma(\xi(1))-\delta\end{aligned}$$

for $r,s\in[0,T]$ and $n\in\{2,\dots,N\}$, for some $\delta>0$. By (5.3), this ensures that

$$\begin{aligned}&\sum_{n=2}^{N-1}\frac{\mathbf{P}\{\xi(1)+\xi(n-s/u)-\xi(\lfloor n-s/u\rfloor)>2u\}}{\mathbf{P}\{\xi(1)>u\}}\\&\qquad\le\sum_{n=2}^{N}\frac{\mathbf{P}\{S_\alpha(2\sigma(\xi(1))-\delta,-1,0)>2u\}}{\mathbf{P}\{\xi(1)>u\}}\\&\qquad\to 0\qquad\text{as } u\to\infty.\end{aligned}$$

Further, by routine estimates, the above scale is at most

$$\begin{aligned}&\|((1+\cdot)^+)^{H-1/\alpha}-((\cdot)^+)^{H-1/\alpha}\\&\quad+((n-s/u+\cdot)^+)^{H-1/\alpha}-((\lfloor n-s/u\rfloor+\cdot)^+)^{H-1/\alpha}\|_\alpha\\&\qquad\le 2^{1/\alpha}[\sigma(\xi(1))+Kn^{H-1}]<2^{1/\alpha}\sigma(\xi(1))\qquad\text{for } s\in[0,T]\text{ and } n\ge N,\end{aligned}$$

for $u$ large enough and some constants $K>0$ and $N\in\mathbb{N}$. Splitting the range of summation $n\in\{N,\dots,\lfloor h/\mathbf{P}\{\xi(1)>u\}\rfloor\}$ into two regions at $n=\lfloor v/u\rfloor^{1-H}$, say, it follows readily, using (5.3), that also the remaining part of (5.9) to treat

$$\lim_{h\downarrow 0}\limsup_{u\to\infty}\sup_{s\in[0,T)}\sum_{n=N}^{\lfloor h/\mathbf{P}\{\xi(1)>u\}\rfloor}\frac{\mathbf{P}\{\xi(1)+\xi(n-s/u)-\xi(\lfloor n-s/u\rfloor)>2u\}}{\mathbf{P}\{\xi(1)>u\}}=0.$$

□

Department of Mathematics
Chalmers University of Technology
412 96 Göthenburg
Sweden
e-mail: palbin@math.chalmers.se